\documentclass[12pt,psamsfonts]{amsart}
\usepackage{epsfig}
\usepackage{graphicx}

\newtheorem{theorem}{Theorem}

\RequirePackage{amsfonts}

\RequirePackage{amssymb}
\RequirePackage{latexsym}
\usepackage[colorlinks=true, hyperfigures=false,
           urlcolor=blue, 
           pdfhighlight =/O]{hyperref}

\makeatletter

\@addtoreset{equation}{section}
\makeatother

\newcommand{\eref}[1]{(\ref{#1})}

\begin{document}
\thispagestyle{empty} \setcounter{page}{1}

\title[Canard Cycles in Global Dynamics]
{Canard Cycles in Global Dynamics}

\author[Alexandre Vidal and Jean--Pierre Fran\c{c}oise]
{Alexandre Vidal and Jean--Pierre Fran\c{c}oise}

\address{ Alexandre Vidal -- Laboratoire Analyse et Probabilit\'es, \newline
\indent Universit\'e d'Evry Val d'Essonne, Evry, France. \newline
\indent \href{mailto:alexandre.vidal@univ-evry.fr}{alexandre.vidal@univ-evry.fr} \newline
\qquad  \newline
\indent Jean--Pierre Fran\c{c}oise -- Laboratoire J.-L. Lions, UMR 7598, CNRS, \newline
\indent Universit\'e P.-M. Curie, Paris6, Paris, France. \newline
\indent \href{mailto:Jean-Pierre.Francoise@upmc.fr}{Jean-Pierre.Francoise@upmc.fr}
}

\thanks{This research is supported by the {\it Agence Nationale de la Recherche} with the project ANAR}

\subjclass{Primary 34C29, 34C25,58F22}

\keywords{Fast-Slow Systems, Canards, Delay to Bifurcation, Transcritical Bifurcation, Relaxation Oscillations}
\date{}
\dedicatory{}

\maketitle

\maketitle

\begin{center}
Submitted for publication on November 15, 2009.
\end{center}

\begin{abstract}
Fast-slow systems are studied usually by ``geometrical dissection" \cite{BR}. The fast dynamics exhibit attractors which may bifurcate under the influence of the slow dynamics which is seen as a parameter of the fast dynamics. A generic solution comes close to a connected component of the stable invariant sets of the fast dynamics. As the slow dynamics evolves, this attractor may lose its stability and the solution eventually reaches quickly another connected component of attractors of the fast dynamics and the process may repeat. This scenario explains quite well relaxation oscillations and more complicated oscillations like bursting. More recently, in relation both with theory of dynamical systems \cite{DR} and with applications to physiology \cite{Desroches, TTFB}, a new interest has emerged in canard cycles. These orbits share the property that they remain for a while close to an unstable invariant set (either singular set or periodic orbits of the fast dynamics). Although canards were first discovered when the transition points are folds, in this article, we focus on the case where one or several transition points or ``jumps" are instead transcritical. We present several new surprising effects like the ``amplification of canards" or the ``exceptionally fast recovery" on both (1+1)-systems and (2+1)-systems associated with tritrophic food chain dynamics. Finally, we also mention their possible relevance to the notion of resilience which has been coined out in ecology \cite{H, LWH, Martin}.
\end{abstract}

\section*{Introduction}

Systems are often complex because their evolution involves different time scales. Purpose of this article is to present several phenomena which can be observed numerically and analyzed mathematically via bifurcation theory.

A first approximation for the time evolution of fast-slow dynamics is often seen as follows. A generic orbit quickly reaches the vicinity of an attractive invariant set of the fast dynamics. It evolves then slowly close to this attractive part until, under the influence of the slow dynamics, this attractive part bifurcates into a repulsive one. Then, the generic orbit quickly reaches the vicinity of another attractive invariant set until it also loses its stability. This approach is a quite meaningful approximation because it explains many phenomena like hysteresis cycles, relaxation oscillations, bursting oscillations \cite{Francoise, Vidal06, Vidal07} and more complicated alternation of pulsatile and surge patterns of coupled GnRH neurons \cite{CF, CV}.

Discovered by E. Benoit, J.-L. Callot, M. and F. Diener (see \cite{BCDD}), canards were first observed in the van der Pol system:
\begin{equation}
\begin{array}{rcl}
\varepsilon \dot{x}&=&y-f(x)=y-(\frac{x^3}{3}+x^2) \\
\dot{y}&=&x-c(\varepsilon)
\end{array}
\label{e1}
\end{equation}
where $c(\varepsilon)$ ranges between some bounds:
\begin{equation}
c_0+\exp \left( -\frac{\alpha}{\varepsilon ^2}\right)<c(\varepsilon )<c_0+ \left( -\frac{\beta}{\varepsilon ^2}\right)
\end{equation}
A canard is an orbit which remains for a while in a small neighborhood of a repulsive branch of the critical manifold $y=f(x)$ i.e. a connected set of repulsive points for the fast dynamics. In the following, we consider the canard phenomenon in its broader sense of delay to the bifurcation of the underlying fast dynamics under the influence of the slow dynamics. More recently, F. Dumortier and R. Roussarie \cite{DR} contributed to the analysis of such orbits by blowing-up techniques and introduced the notion of canard cycle. There are now several evidences showing the relevance of this notion to explain experimental facts observed in physiology (see \cite{Desroches, TTFB}).

The systems presented here display anomalous long delay to ejection from the repulsive part of the fast dynamics.

\section{Enhanced delay and canard cycles of planar (1+1)--dynamics}

\subsection{Dynamical Transcritical Bifurcation} \label{DynTrans}
\qquad \newline

The classical transcritical bifurcation occurs when the parameter $\lambda$
in the equation:
\begin{equation}
\dot{x}=-{\lambda}x+x^2 \label{e3}
\end{equation}
crosses $\lambda=0$. Equation \eref{e3} displays two equilibria, $x=0$
and $x=\lambda$. For $\lambda>0$, $x=0$ is stable and $x=\lambda$
is unstable. After the bifurcation, $\lambda<0$, $x=0$ is stable
and $x=\lambda$ is unstable. The two axis have ``exchanged" their
stability.

The terminology ``Dynamical Bifurcation'' (due to R. Thom) refers to the situation
where the bifurcation parameter is replaced by a slowly varying
variable. In the case of the transcritical bifurcation, this
yields:
\begin{equation}
\begin{array}{rcl}
\dot{x}&=&-yx+x^2 \\
\dot{y}&=&-\varepsilon
\end{array}
\label{e4}
\end{equation}
where $\varepsilon$ is assumed to be small. This yields:
\begin{equation*}
\dot{x}=-(-\varepsilon t+y_0)x+x^2, \qquad (y_0=y(0))
\end{equation*}
which is an integrable equation of Bernoulli type. Its solution is:
\begin{equation*}
x(t)=\frac{x_0\exp [-Y(t)]}{1-x_0\int_0^t \exp [-Y(u)]du}, \qquad (x_0=x(0)) \\
\end{equation*}
where:
\begin{equation*}
Y(t)=\int_0^t y(s)ds=\int_0^t (-\varepsilon s+y_0)ds=-\varepsilon \frac{t^2}{2}+y_0 t
\end{equation*}

\begin{figure}[htb]
\centering
\includegraphics[width=8cm]{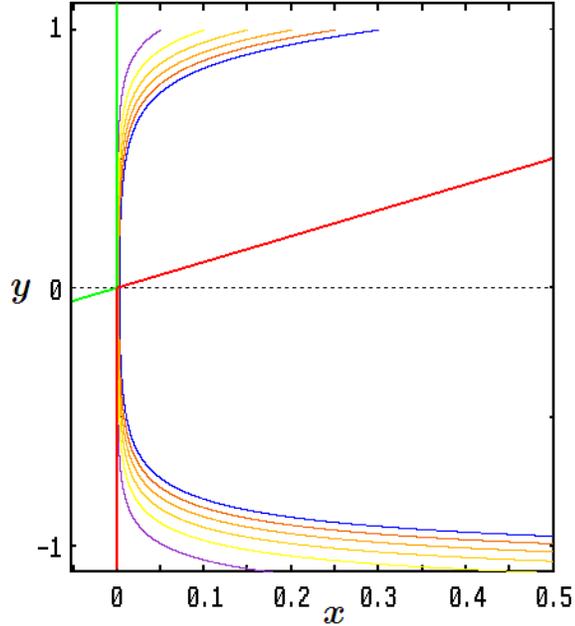}
\caption{Orbits of system \eref{e4} starting from $(0.05 i, 1)$, for $i=1,..,6$, with $\varepsilon =0.1$. The critical set -- formed by the singular points of the fast dynamics -- is shown in green and red. On this set, the red points are repulsive for the fast dynamics and the green points are attractive. Each orbit reaches a neighborhood of the attractive manifold $x=0, y>0$, goes down along $x=0$. When $y$ becomes negative, despite the repulsiveness of $x=0, y<0$, each orbit remains for a very long time close to $x=0$ before drifting away.}
\label{Delay}
\end{figure}

If we fix an initial data $(x_0,y_0)$, $y_0>0$, $0<x_0<y_0/2$, and we consider
the solution starting from this initial data, we find easily that it takes
time $t=y_0/\varepsilon $ to reach the axis $y=0$. If $x_0$ is quite
small, that means the orbit stays closer and closer of the
attractive part of the critical manifold until it reaches the axis $x=y$ and then
coordinate $x$ starts increasing. But now consider time
$cy_0/\varepsilon $, $1\leq c\leq 2$. Then, a straightforward computation shows that:
\begin{equation*}
Y(t)=c\left(1-\frac{c}{2}\right)\frac{y_0^2}{\varepsilon} =\frac{k}{\varepsilon }
\end{equation*}
and:
\begin{equation*}
x(t)=O\left(\frac{x_0{\rm e}^{-\frac{k}{\varepsilon}}}{1-\frac{2x_0}{y_0}}\right)
\end{equation*}

This shows that, despite the repulsiveness of the half-line $x=0, y<0$ for the fast dynamics,
the orbit of \eref{e4} remains for a very long time close to $x=0$, indeed
$x(cy_0/\varepsilon )<2x_0$ for $\varepsilon $ small enough (see Figure \ref{Delay}). Note that, afterwards for larger time, the orbit blows away
from this repulsive axis. This phenomenon, although quite simply
explained, is of the same nature as the delay to bifurcation
discovered for the dynamical Hopf bifurcation, see for instance \cite{ErnBaRinz, CDD, Ern, Z}. This well-known effect is instrumental in the systems we study in this article. Some related work has been done in computing entry-exit relation for the passage near single turning points (see \cite{Ben, DeDu}).

\subsection{Enhanced delay to bifurcation} \label{SectDTC}
\qquad \newline

Consider the system:
\begin{equation}
\begin{array}{rcl}
\dot{x}&=&(1-x^2)(x-y) \\
\dot{y}&=&\varepsilon x
\end{array}
\label{DTC}
\end{equation}
The critical set -- defined as the set of singular points of the fast dynamics -- is the union of the straight lines $x=-1$, $x=1$, and $y=x$. In the following, we call ``slow manifolds'' the two straight lines $x=-1$ and $x=1$, as they are invariant for the critical system $\{(1-x^2)(x-y)=0, \dot{y} =x\}$. A quick analysis shows that, as the slow variable $y$ (considered as a parameter) evolves, the fast system undergoes two transcritical bifurcations: at $x=-1$ for $y=-1$, at $x=1$ for $y=1$.

As we recalled in Subsection \ref{DynTrans}, a typical orbit starting from an initial data $(x_0,y_0)$, $\left\vert x_0 \right\vert <1$, close to $(x=-1, y>-1)$, first goes down along $x=-1$ and displays a ``delay'' along the repulsive part of the slow manifold $(x=-1,y<-1)$. Then, under the influence of the fast dynamics, it quickly reaches the attractive part $(x=1,y<1)$ and moves upward to the other transcritical bifurcation point. There, it again displays a delay along the repulsive part $(x=1,y>1)$. Then, it quickly reaches $(x=1, y>-1)$ and starts again. Hence, there is a mechanism of successive enhancements of the delay to bifurcation after each oscillation generated by the hysteresis. We proved in \cite{FPV} the:
\begin{theorem} \label{ThNoLCDTC}
For all initial data inside the strip $-1<x<1$, for all $\delta$
and for all $T$, the corresponding orbit spends a time larger than
$T$ within a distance less than $\delta$ to the repulsive part of
the slow manifolds.
\end{theorem}

We also proved in \cite{FPV} the:
\begin{theorem} \label{ThAsymptDTC}
Given any initial data $(x_0,y_0)$ outside the strip $\left\vert  x\right\vert \leq
1$, the corresponding orbit is asymptotic to $y=x$.
\end{theorem}

\begin{figure}[htb]
\centering
\includegraphics[width=14cm]{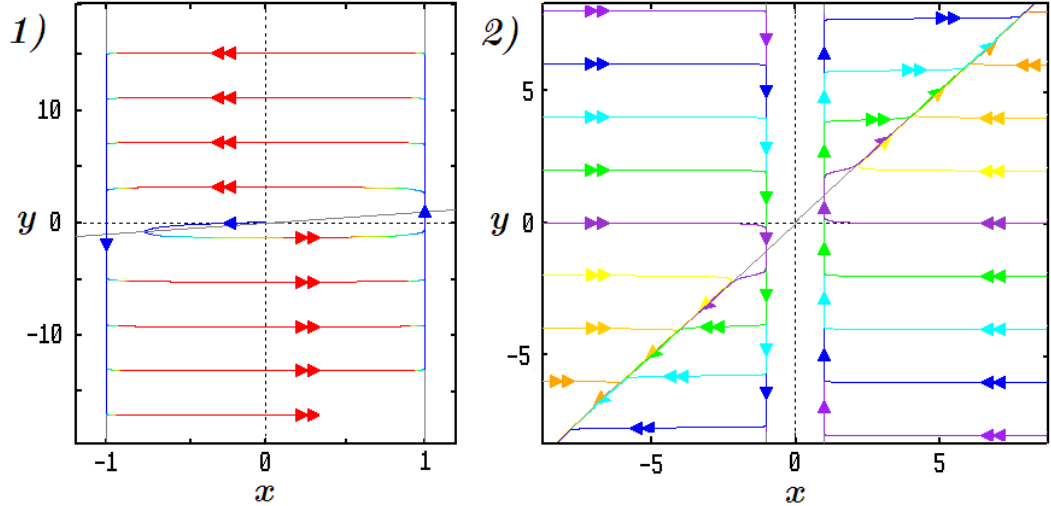}
\caption{Typical orbits of \eref{DTC} with $\varepsilon =0.1$ starting from a point near the origin in panel 1) and from various initial datas outside the strip $\left\vert x \right\vert <1$ in panel 2). Double arrows are added on the fast parts of the orbits and single arrows on the slow parts. \newline
Panel 1): The delay to the transcritical bifurcation undergone by the orbit is enhanced at each passage. Consequently, the time taken to escape a given neighborhood of $x=-1, y<-1$ on one hand and a given neighborhood of $x=1, y>1$ on the other hand is longer for each oscillation than for the preceding one, although these half-lines are repulsive for the fast dynamics. \newline
Panel 2): Each orbit is asymptotic to $y=x$. \newline
These numerical simulations were performed with order 4 Runge-Kutta integration scheme (absolute error value: $10^{-9}$ ; relative error value: $10^{-12}$ ; mean integration step: $10^{-5}$).}
\label{DTCSimul}
\end{figure}

\subsection{Structural stability of the enhanced delay} \qquad \newline

The theory of the structural stability of fast-slow systems remains to be found. In this subsection, we actually adopt a very pragmatic approach and restrict ourselves to numerical simulations. We perturb the system by adding an arbitrary small perturbation generated by a chaotic system, more precisely the R\"{o}ssler system \cite{Rossler}.
\begin{equation} \label{DTCP}
\begin{array}{rcl}
\varepsilon \dot{x}&=&(1-x^2)(x-y)+{\alpha}n(t) \\
\dot{y}&=&x
\end{array}
\end{equation}
where $n(t)$ is the first variable of a R\"{o}ssler system:
\begin{equation}
\begin{array}{rcl}
\dot{n} &=& -u-v \\
\dot{u} &=& n+au \\
\dot{v} &=& b+(n-c)v
\end{array}
\end{equation}
with $a=0.1$, $b=0.1$, and $c=14$
with some fixed initial data.

We choose $\alpha$ quite small. The perturbation does not change the behavior of the typical orbit for a while. It starts to oscillate between $x=1$ and $x=-1$: it remains for a long time alternatively near each of these lines and undergoes fast motions after ejection. But, as the orbit approaches closely these axes (within the distance of
$\alpha$), the small chaotic perturbation becomes operating and moves slightly the orbit outside of the strip. After crossing one of the straight lines $x=-1$ or $x=1$, the orbit moves quickly to the axis $y=x$ as previously shown.

\begin{figure}[htb]
\centering
\includegraphics[width=8.5cm]{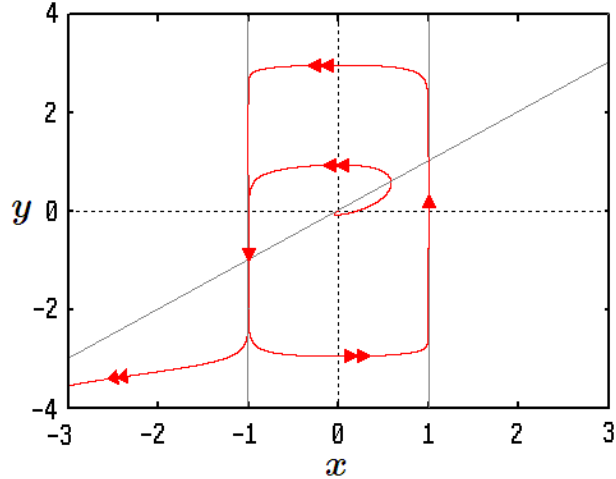}
\caption{Example of an orbit of \eref{DTCP}, with $\varepsilon =0.5$ starting from $(x,y)=(10^{-3},0)$ for $\alpha =10^{-4}$. Double arrows are added on the fast parts of the orbit and single arrows on the slow parts.}
\label{DTCperturb}
\end{figure}

Figure \ref{DTCperturb} displays an example of such an orbit starting from $(0.001,0)$ for $\alpha =10^{-4}$. Despite the very small magnitude of the perturbation -- less than $3.10^{-3}$ -- the orbit exits the strip $\left\vert x \right\vert \leq 1$ before the third oscillation around the origin. This simulation shows the loss of the oscillatory behavior as, under the influence of a perturbation of very small magnitude, the orbit crosses one of the slow manifolds as soon as $x$ is sufficiently close to $-1$ or $1$. We have thus obtained numerical evidence that the system is not structurally stable. 

It is interesting to point out the fact that the numerical simulation of \eref{DTC} highlights also this structural unstability. In fact, the integration scheme itself -- whatever the integration step and the tolerances -- provides errors. They give rise in the long run to intrinsic approximation: after many oscillations, the orbit starting from $\left\vert x \right\vert <1$ passes so close to the slow manifolds that the numerical integration leads, sooner or later, to the approximation $x=1$ or $x=-1$ or can eventually cross one of these lines.

\subsection{Canard cycle for a double transcritical system} \label{SectDTCBB}
\qquad \newline

In this subsection, we introduce a new system inspired by the preceding \eref{DTC}:
\begin{equation}
\begin{array}{rcl}
\dot{x}&=&(1-x^2)(x-y) \\
\dot{y}&=&\varepsilon x(y+b)(a-y)
\end{array}
\label{DTCBB}
\end{equation}
where $a,b>1$ are parameters.
The critical set is the same as the critical set of \eref{DTC}, formed by the three straight lines of equation $x=-1$, $x=1$ and $y=x$. The dynamics displays again two transcritical bifurcations: at $x=-1$ for $y=-1$, at $x=1$ for $y=1$. But, the novelty is in the factors $(a-y)$ and $(y+b)$ in $\dot{y}$ which yield bounded orbits. We restrict the phase space to the compact set:
\[
K=\{(x,y)|-1 \leq x \leq 1, -b \leq y \leq a \}
\]

As for \eref{DTC}, the origin is a repulsive focus of \eref{DTCBB}. Except the origin, the singular points of \eref{DTCBB} in $K$ are the summits of the rectangular boundary of $K$, $(-1,-b)$, $(1,-b)$, $(1,a)$, $(-1,a)$, and are all of saddle type. As the straight lines $x=-1$, $x=1$, $y=a$, $y=-b$ are invariant under the flow of \eref{DTCBB}, we deduce the stable and unstable manifolds of each saddle:
\vspace{0.5cm}
\begin{center}
\begin{tabular}{|l|ll|ll|}
\hline
saddle & stable & manifold & unstable & manifold\\
\hline
$(-1,-b)$ & $x=-1$, & $y<a$  &  $y=-b$, & $x<1$   \\
\hline
$(1,-b)$  & $y=-b$, & $x>-1$  &  $x=1$, & $y<a$    \\
\hline
$(1,a)$   & $x=1$, & $y>-b$  &  $y=a$, & $x>-1$    \\
\hline
$(-1,a)$  & $y=a$, & $x<1$   &  $x=1$, & $y>-b$   \\
\hline
\end{tabular}
\label{SaddlesDTCBB}
\end{center}
\vspace{0.5cm}
Consequently, the interior of $K$ is invariant under the flow and contains only a repulsive focus. The boundary of $K$ is a graphic ($\omega $-limit set formed by the union of the saddles with their separatrices).

Following the study on the dynamical transcritical bifurcation, one expects that the value of $y$ along an orbit oscillates between $-b$ and $\min (a,b+2)$ if $a \geq b$ or between $-\min(b, a+2)$ and $a$ if $a < b$. However, the exponential attraction of $x=-1$, $-1<y<a$ (resp. $x=1$, $-b<y<1$), formed by attractive points of the fast dynamics, produces a delay to bifurcation and keeps the orbit near $x=-1$, $-b<y<-1$ (resp. $x=1$, $1<y<a$), formed by repulsive points of the fast dynamics. The delay could be great enough so that the orbit reaches a small neighborhood of $\dot{y}=0$. Hence, after the fast motion, the orbit tracks the other manifold $x=1$, $-b<y<1$ (resp. $x=-1$, $-1<y<a$) even closer than during the preceding passage. The delay obtained afterwards is then again enhanced, as the orbit approaches even closer the slow manifolds. Hence, an orbit starting near the origin, after several oscillations, reaches alternatively a small neighborhood of both $y=a$ and $y=-b$, whatever the values of $a$ and $b$.

\begin{center}
\begin{figure}[htb]
\begin{tabular}{cc}
\includegraphics[width=0.47\textwidth]{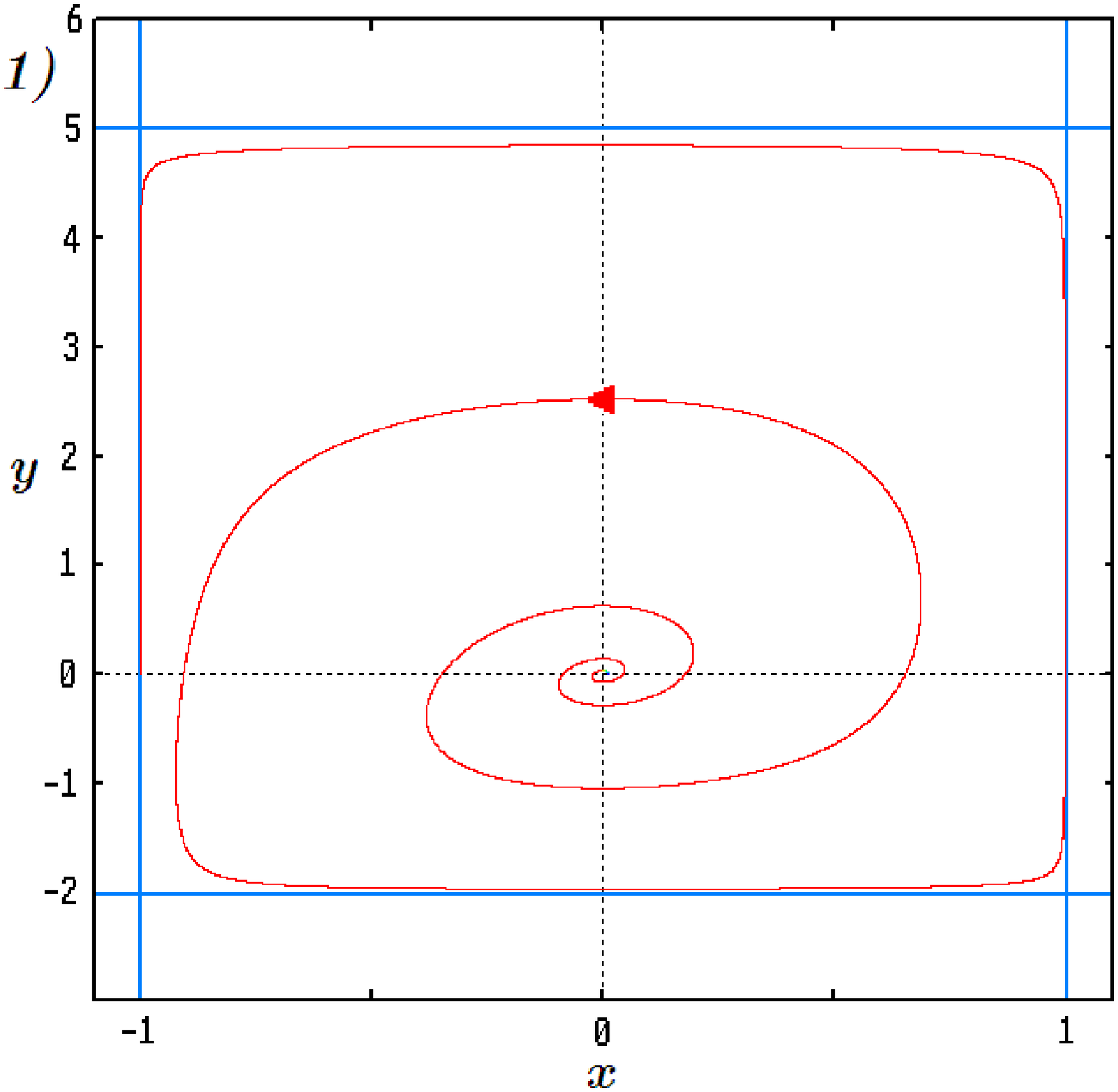} & \includegraphics[width=0.47\textwidth]{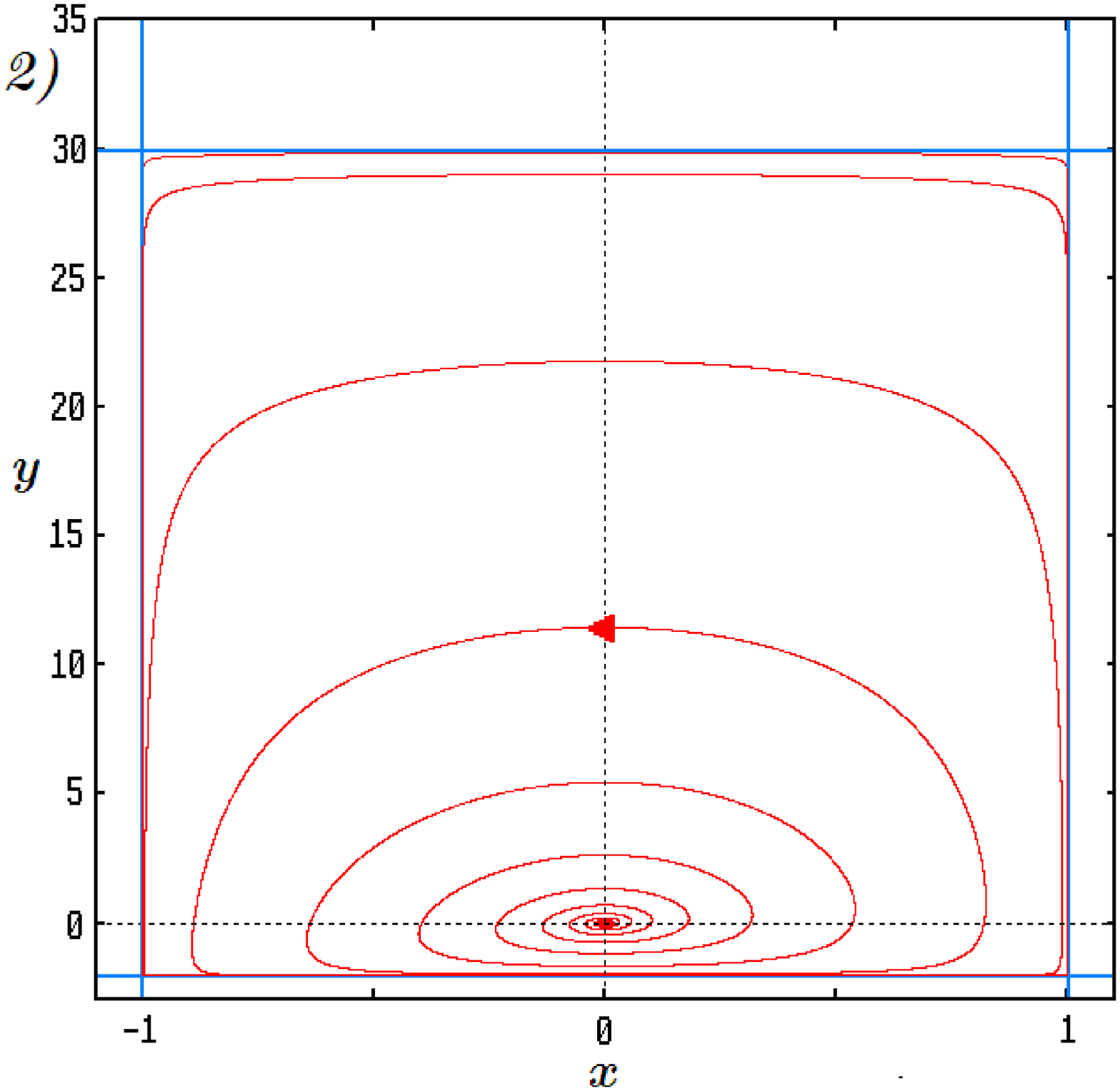}
\end{tabular}
\caption{Orbits (red) of \eref{DTCBB} for $\varepsilon =0.5$, $b=2$ and respectively $a=5$ in panel 1) and $a=30$ in panel 2). In both cases, the initial data is $(10^{-3},0)$. Comparison of panels 1) and 2) shows that, whatever the values of $a$ and $b$, the orbit enters alternatively small neighborhood of $x=-1$, $x=1$, $y=-b$ and $y=a$ (blue straight lines). These numerical simulations were performed with order 4 Runge-Kutta integration scheme (absolute error value: $10^{-9}$ ; relative error value: $10^{-12}$ ; mean integration step: $10^{-6}$).}
\label{DTCBBSimul}
\end{figure}
\end{center}

Figure \ref{DTCBBSimul} illustrates this global behavior with fine step numerical simulations for various sets of parameter values. As expected, the smaller $\varepsilon $, the faster the orbit reaches a given neighborhood of $x=-1$ (resp. $x=1$). Starting near the origin, if $\varepsilon $ is small (for instance, equal to $0.5$ for the simulations presented in Figure \ref{DTCBBSimul}), a few oscillations suffice to obtain a part close to the graphic.

Thus, the orbits of system \eref{DTCBB} generate a new type of canards. On the contrary of the canards discovered in the solutions of the van der Pol system, all the orbits of system \eref{DTCBB} starting from $\mathring{K}\backslash \{(0,0)\}$ display delays to bifurcation. Moreover, we have a very simple way to modulate these delays by choosing the parameter values $a$ and $b$. It is worth noticing that, as in system \eref{DTC}, small perturbations may provoke dramatic changes in the orbits.

Now, we explore the ejection mechanism -- leading to this type of oscillations -- using the conjugacy of \eref{DTCBB} with an appropriate semi-local form and the analysis of the transition near the saddle $(x,y)=(-1,-b)$ (the case of $(x,y)=(1,a)$ is similar). First, we translate the origin to the saddle via the change of variables $\{X=x+1, Y=y+b\}$ to obtain the new form of the system:
\begin{equation}
\begin{array}{rcl}
\dot{X}&=&X(X-Y+b-1)(2-X) \\
\dot{Y}&=&\varepsilon Y(X-1)(a+b-Y)
\end{array}
\label{DTCBBTr}
\end{equation}
Note that the left point of transcritical bifurcation of the fast dynamics $(x,y)=(-1,-1)$ reads now $(X,Y)=(0,b-1)$.

Consider an initial data $(\bar{X}, \bar{Y})$ with $\bar{X} \in]0,1], b-1<\bar{Y}<a+b$. As explained above, for $\varepsilon $ small, the corresponding orbit of \eref{DTCBBTr} reaches quickly the vicinity of $X=0, y>b-1$, $X$ decreases very quickly while $Y$ remains almost constant. Hence, along this fast branch of the dynamics, variable $Y$ is given by $\bar{Y}+O(\varepsilon )$ until $X=O(\varepsilon)$. Afterwards,
approximation of the dynamics near $X=0$ by the slow motion and direct integration yield an approximation for the time needed to reach $Y=b-1$:
\[
T_{\bar{Y} \rightarrow b-1} \underset{\varepsilon \rightarrow 0}{=} \frac{1}{\varepsilon (a+b)} \ln \left[ \frac{(a+1)\bar{Y}}{(b-1)(a+b-\bar{Y})} \right]+O(1)
\]
The leading term of the $x$-component $C \exp (-k/\varepsilon)$ displays:
\begin{eqnarray}
C &=& (b-1)^{2\frac{b-1}{a+b}} (a+1)^{2\frac{a+1}{a+b}} \label{paramC}\\
k &=& \frac{2}{a+b} \left[ (b-1) \ln (\bar{Y}) + (a+1) \ln (a+b-\bar{Y}) \right] >0 \label{paramk}
\end{eqnarray}
We now study the delay to bifurcation and consider the entry of the orbits coming from above $Y=b-1$. Hence, we consider initial values of $X$ which are exponentially small with respect to $\varepsilon $: $C \exp (-k/\varepsilon)$.

\begin{figure}[htb]
\centering
\includegraphics[height=7cm]{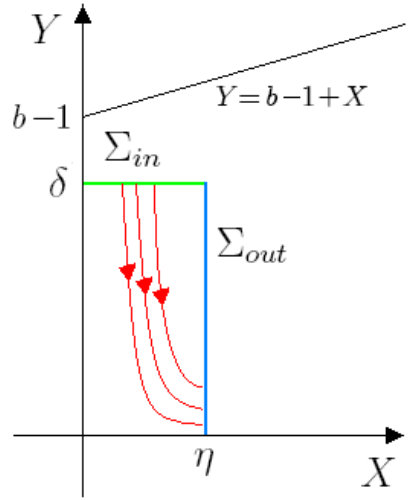}
\caption{The transition function induced by the flow of \eref{DTCBBTr} is well-defined from $\Sigma _{in}$ into $\Sigma _{out}$.}
\label{LocAn}
\end{figure}

To specify the transition induced by the flow near the saddle, we consider the two sections:
\begin{eqnarray}
\Sigma_{in} &=& \left\{(X,\delta )|0<X \leq \eta \right\} \\
\Sigma_{out} &=& \left\{(\eta,Y)|0<Y \leq \delta \right\}
\end{eqnarray}
with $\eta >0$ a small fixed parameter and $0<\delta \leq b-1$. We note:
\[
U=\left\{(X,Y)|0<X \leq \eta, 0<Y \leq \delta \right\}
\]
the rectangle delimited by $\Sigma_{in}$, $\Sigma_{out}$ and the stable and unstable manifolds of the saddle (see Figure \ref{LocAn}).

Hence, as $\dot{Y}<0$ in $\{(X,Y)|0<X<1,0<Y<X+b-1 \}$ and $\delta \leq b-1$, any orbit starting from $\Sigma _{in}$ enters $U$ and escapes from $U$ through $\Sigma _{out}$. Thus, the transition function induced by the flow is well-defined from $\Sigma _{in}$ into $\Sigma_{out}$.

As $X$ is small in $U$ (smaller than $\eta $), system \eref{DTCBBTr} is conjugated to:
\begin{equation}
\begin{array}{rcl}
\dot{X}&=&2 (b-1) X-2XY \\
\dot{Y}&=&\varepsilon Y
\end{array}
\label{DTCBNL}
\end{equation}
Note that $(X,Y)=(0,b-1)$ is also a point of transcritical bifurcation for this system.
Direct integration provides the orbit from $(X_0, \delta )$:
\begin{eqnarray}
X(t) &=& X_{0}\exp \left[ 2(b-1)t+\frac{2\delta }{\varepsilon }\left({\rm e}^{-\varepsilon t}-1\right) \right] \nonumber\\
Y(t) &=& \delta {\rm e}^{-\varepsilon t} \nonumber
\end{eqnarray}

As explained previously, we consider initial data of type $X_0 = C \exp (-k/\varepsilon)$, where $k, C>0$. If $\eta <C$, for $\varepsilon $ small enough, i.e.:
\begin{equation*}
\varepsilon < - \frac{k}{\ln \frac{\eta}{C}}
\end{equation*}
$(X_0,\delta )$ lies in $\Sigma _{in}$. If $\eta \geq C$, as $\varepsilon ,k>0$, all values of $\varepsilon $ fulfill this property.

The time $T$ needed to go from $(X_0,\delta)$ to $(\eta , Y(T)) \in \Sigma_{out}$ along the flow is the solution of:
\begin{equation*}
2(b-1)T+\frac{2\delta }{\varepsilon }{\rm e}^{-\varepsilon T}=\frac{k}{\varepsilon } + \ln \frac{\eta}{C}+\frac{2\delta}{\varepsilon}
\end{equation*}
The solution $T$ can be expressed via the ``Lambert function'' $\mathcal{W_L}$ -- inverse function of $w\rightarrow w {\rm e}^{w}$ (see \cite{CGHJK}). This yields:
\begin{multline}
T=\frac{1}{\varepsilon } \mathcal{W_L} \left( -\frac{\delta}{b-1} \exp \left( -\frac{\epsilon \ln(\eta /C)+k+2\delta }{2(b-1)} \right) \right) \\
     + \frac{1}{2(b-1)} \left( \ln \frac{\eta }{C} + \frac{k+2\delta }{\varepsilon} \right)
\label{TempsTrans}
\end{multline}
As expected, the transition time is $O(1/\varepsilon)$ and the leading term of $T$ is:
\begin{equation}
T \underset{\varepsilon \rightarrow 0}{=} \frac{1}{\varepsilon} \left[ \mathcal{W_L} \left( -\frac{\delta}{b-1} {\rm e}^{-\frac{2\delta +k}{2(b-1)}} \right) + \frac{2 \delta+k}{2(b-1)} \right] + O(1)
\label{TTrans}
\end{equation}

This displays the transition function $(X_0, \delta) \rightarrow (\eta, Y_{out}(X_0))$ with:
\begin{multline}
Y_{out}(C {\rm e}^{-\frac{k}{\varepsilon }}) = \\
   \delta \exp \left[ -\mathcal{W_L} \left( -\frac{\delta}{b-1} \exp \left( -\frac{\epsilon \ln(\eta /C)+k+2\delta }{2(b-1)} \right) \right) \right.\\
     \left. + \frac{1}{b-1} \left( \ln \frac{\eta }{C} + \frac{k+2\delta }{\varepsilon} \right) \right]
\end{multline}
This shows that, even if the transition time tends to $+\infty $, the exit function displays as $O(1)$-leading term:
\begin{multline}
Y_{out}(C {\rm e}^{-\frac{k}{\varepsilon }}) \underset{\varepsilon \rightarrow 0}{=} Y_{out}^{0}(C {\rm e}^{-\frac{k}{\varepsilon }}) + O(\varepsilon) \\
     \underset{\varepsilon \rightarrow 0}{=} \delta \exp \left[ -\mathcal{W_L} \left( -\frac{\delta}{b-1} {\rm e}^{-\frac{2\delta +k}{2(b-1)}} \right) -\frac{2\delta +k}{2(b-1)} \right] +O(\varepsilon)
\label{DLYout}
\end{multline}

Hence, the longer the orbit has remained near the slow manifold, the stronger the ejection (to reach the neighborhood of $(\eta, Y_{out}^{0})$ given by \eref{DLYout}) is. This property of the orbits -- staying a very long time near the repulsive part of the slow manifold without squashing on the $\dot{y}$-nullcline -- together with the increasing strength of ejection is what we call the ``exceptionally fast recovery''. Similar study can be done for the other transition of \eref{DTCBB} and this leads to similar results (where $a$ and $b$ are exchanged).

Finally, note that, for given values of parameter $a$, $b$, $\varepsilon $, the values of parameters $k$ and $C$ in \eref{TTrans} and \eref{DLYout} are approximated for the global dynamics using \eref{paramC} and \eref{paramk}.

\section{Saddle-node transcritical ejection in a prey-predator-superpredator model}

\subsection{Tritrophic food chain dynamics} \qquad \newline

In the late seventies, interest in the mathematics of tritrophic food chain models (composed of prey, predator and superpredator) appear (see, for instance, \cite{FW, Ga}). Related predator-prey models with parasitic infections were studied later \cite{HF}. In the nineties, in \cite{HP, MR2} and \cite{KH}, the existence of chaotic attractors was discussed. There are many more recent contributions, that we can not refer in more details (see, for instance, \cite{KBK, Vidal07}).

We investigate here the following:
\begin{eqnarray}
\frac{dU}{dT}& = &U\left( R\left( 1-\frac{U}{K}\right) -\frac{A_{1}V}{B_{1}+U}\right)  \nonumber \\
\frac{dV}{dT}& = &V\left( E_{1}\frac{A_{1}U}{B_{1}+U}-D_{1}-\frac{A_{2}W}{B_{2}+V}\right)  \label{sys1} \\
\frac{dW}{dT}& = &\varepsilon W\left( E_{2}\frac{A_{2}V}{B_{2}+V}-D_{2}\right)  \nonumber
\end{eqnarray}
which represents the interactions between three populations $U$, $V$ and $W$. The
variable $U$ stands for the prey, $V$ for its predator and $W$ for a superpredator of $V$.
The threshold constant $K>0$ and the intrinsic growth rate of the prey $R>0$
characterize the logistic evolution of $U$.

The predator--prey interactions are described by two Holling type II
factors defined by the positive parameters:
\begin{align*}
A_{j}& :\text{the maximum predation rates} \\
B_{j}& :\text{the half-saturation constants} \\
D_{j}& :\text{the death rates} \\
E_{j}& :\text{the efficiencies of predation}
\end{align*}
$j=1$ relates to the predator $V$, $j=2$ relates to the superpredator $W$.
It is assumed that the evolution of the superpredator is slower
than those of the predator and the prey. Then we introduce different time scales by means of the
constant $0< \varepsilon \ll 1$. 

\subsection{Bifurcations of the fast dynamics} \qquad \newline

The global behavior of this system (existence of global periodic orbit, bifurcations of limit cycles, early crisis in the predator membership) has been studied in \cite{Vidal06} and \cite{Vidal07}. We recall briefly the classical situation of interest for us for which the system is bistable.

In order to obtain a simpler and more useful analytic form, Klebanoff and Hastings\- proposed in \cite{KH} the following rescalings:
\begin{equation}
x=\frac{U}{K},\qquad y=\frac{V}{KE_{1}},\qquad z=\frac{W}{KE_{1}E_{2}}
,\qquad t=RT  \label{Rescal}
\end{equation}
which yields:
\begin{eqnarray}
\dot{x} &=&x\left( 1-x-\frac{a_{1}y}{1+b_{1}x}\right) = f(x,y,z)  \notag \\
\dot{y} &=&y\left( \frac{a_{1}x}{1+b_{1}x}-d_{1}-\frac{a_{2}z}{1+b_{2}y} \right) = g(x,y,z)  \label{sys2} \\
\dot{z} &=&\varepsilon z\left( \frac{a_{2}y}{1+b_{2}y}-d_{2}\right) = h(x,y,z)  \notag
\end{eqnarray}
where $a_j, b_j, d_j, j=1,2$ are positive parameters (see \cite{KH} or \cite{Vidal06} for their expressions in function of $A_j, B_j, D_j, E_j, R$ and $K$.
All axes and faces of the positive octant $\mathbb{R} _{+}^{3}$ are invariant sets of \eref{sys2}. Thus, we limit the phase space to this positive octant.

Considering the slow variable $z$ as a parameter, we describe the sequence of bifurcations undergone by the two-dimensional fast dynamics. To this purpose, it is convenient to introduce the critical set:
\begin{equation}
C=\{(x,y,z)\in \mathbb{R}_{+}^{3}|f(x,y,z)=0,g(x,y,z)=0\}
\end{equation}
formed by the singular points of the so-called Boundary-Layer System (BLS), obtained from \eref{sys2} with $\varepsilon =0$. For any point $(\tilde{x}, \tilde{y}, \tilde{z}) \in C$, $(\tilde{x}, \tilde{y})$ is a singular point of the fast dynamics for $z=\tilde{z}$. In the following, we note $\pi $ the projection from $\mathbb{R}_{+}^{3}$ into $\mathbb{R}_{+}^{2}$: $\pi (x,y,z) =(x,y)$.

First, we assume:
\begin{equation}
G=a_{1}-d_{1}(1+b_{1})>0
\label{StAss}
\end{equation}
to ensure that the singular point:
\begin{equation}
\left( \frac{d_{1}}{G+d_{1}},\frac{G}{\left( G+d_{1} \right) ^{2}},0\right)
\end{equation}
lays in the phase space $\mathbb{R} _{+}^{3}$.

The critical set writes $C=\Delta \cup \mathcal{L}$ where:
\begin{eqnarray}
\Delta  &=&\{(1,0,z)|z\in \mathbb{R}_{+}\} \\
\mathcal{L} &=&\left\{ \left( x,y_{\mathcal{L}}(x),z_{\mathcal{L}}(x)\right) \in \mathbb{R} _{+}^{3}|x\in [0,1]\right\} 
\end{eqnarray}
and:
\begin{align}
y_{\mathcal{L}}(x)& =\frac{1}{a_{1}}\left( 1-x\right) \left( 1+b_{1}x\right)  \\
z_{\mathcal{L}}(x)& =\frac{\left( a_{1}x-d_{1}\left( 1+b_{1}x\right) \right) \left( a_{1}+b_{2}(1-x)(1+b_{1}x)\right) }{a_{1}a_{2}(1+b_{1}x)}
\end{align}
Note that $\Delta $ always intersects $\mathcal{L}$ at the point $T=(1,0,z_{T})$
where:
\begin{equation}
z_{T}=\frac{G}{a_{2}(1+b_{1})} >0 \label{zT}
\end{equation}

In the following, we assume that $d_{1}$ is small enough so that:
\begin{equation}
\text{there is a point }x_{P}>0\text{ so that }z_{\mathcal{L}}^{\prime }(x_{P})=0.
\label{PaboveT}
\end{equation}
We note $P=(x_P,y_P,z_P)=(x_P,y_{\mathcal{L}}(x_{P}),z_\mathcal{L}(x_P))$. Let us remark that, under assumption \eref{PaboveT}, $z_{T}<z_{P}$.
Consequently, $\mathcal{L}$ is $\cap $-shaped in $\mathbb{R} _{+}^{3}$ (cf. Figure \ref{EnsCrit}) and we note:
\begin{itemize}
\item for each $z\in \left] z_{T},z_{P}\right[ $, $S_z$ the unique point $(x,y_{\mathcal{L}}(x),z)\in \mathcal{L}$ \newline
such that $x> x_P$ ; $\mathcal{L}_S=\underset{z_{T}<z<z_{P}}{\cup }\{S_{z}\}$ ;
\item for each $z\in \left[ 0,z_{P}\right[ $, $R_z$ the unique point $(x,y_{\mathcal{L}}(x),z)\in \mathcal{L}$ \newline
such that $x< x_P$ ; $\mathcal{L}_{\pm}=\underset{0<z<z_{P}}{\cup }\{R_{z}\}$.
\end{itemize}
Thus, $\mathcal{L}={T} \cup \mathcal{L}_S \cup {P} \cup \mathcal{L}_{\pm}$ (see Figure \ref{EnsCrit}).

The points $\pi (S_z)$ are saddle type singular points of the fast dynamics, $(1,0)$ is a saddle for $z<z_T$ and an attractive node for $z>z_T$.
Additionally, we assume that $d_{1}$ is small enough such that there exists a point $z_{H}\in ]z_{T},z_{P}[$ so that:
\begin{eqnarray*}
\text{for }&0\leq z<z_{H},&\pi (R_z)\text{ is a repulsive focus}, \\ 
\text{for }&z_{H}\leq z<z_{P},&\pi (R_z)\text{ is an attractive focus}.
\end{eqnarray*}
Figure \ref{EnsCrit} displays $C$ and its splitting according to the nature of the singular points for the fast dynamics. It can be seen as a bifurcation diagram of the fast dynamics as the bifurcation parameter $z$ varies. Hence, as $z$ decreases, the following sequence of bifurcations occurs:
\begin{figure}[htb]
\centering
\includegraphics[height=9cm]{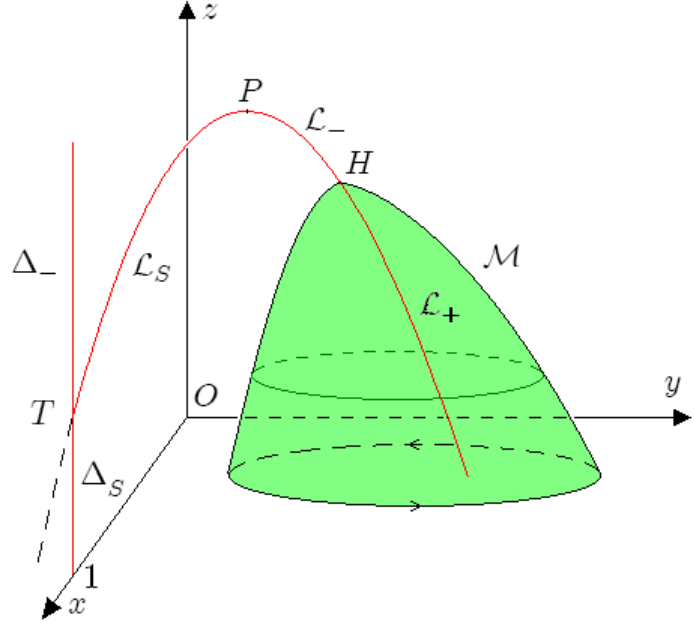}
\caption{Critical set of \eref{sys2} (in red) -- set of singular points of the Boundary-Layer System (BLS) obtained by setting $\varepsilon =0$. It is splitted according to the nature of the singular points for the fast dynamics with the corresponding value of $z$: $\Delta _-$ and $\mathcal{L}_-$ formed by attractive nodes, $\Delta _S$ and $\mathcal{L}_S$ by saddles, $\mathcal{L}_+$ by repulsive foci. The points $P$, $H$, $T$ are respectively the points of saddle-node, Hopf and saddle-node transcritical bifurcation of the fast dynamics. For $\tilde{z} \in [0,z_H[$, the repulsive focus $\pi (R_{\tilde{z}}) \in \mathcal{L}_+$ is surrounded by an attractive limit cycle of the fast dynamics. The union $\mathcal{M}$ of these limit cycles (in green) is an invariant attractive manifold for the (BLS).}
\label{EnsCrit}
\end{figure}
\begin{itemize}
\item for $z>z_P$, the attractive node $(1,0)$ is the unique singular point.
\item as $z=z_P$, an inverse saddle-node bifurcation occurs at $(x_P,y_P)$.
\item for $z_H<z<z_P$, there are three singular points: the attractive node $(1,0)$, the saddle $\pi (S_z)$, the attractive focus $\pi (R_z)$.
\item as $z=z_H$, $\pi (R_z)$ undergoes a supercritical Hopf bifurcation.
\item for $z_T<z<z_H$, there are three singular points: the attractive node $(1,0)$, the saddle $\pi (S_z)$, the repulsive focus $\pi (R_z)$ surrounded by an attractive limit cycle, born from the Hopf bifurcation.
\item as $z=z_T$, a saddle-node transcritical bifurcation occurs at $(1,0)$ ($\pi (S_z)$ and $(0,1)$ exchange their stability).
\item for $0\leq z< z_T$, there are two singular points in the positive octant: the saddle $(0,1)$ and the repulsive focus surrounded by an attractive limit cycle.
\end{itemize}
We note $\mathcal{M}$ the union in $\mathbb{R}_{+}^{3}$ of the planar attractive limit cycles surrounding the points $R_z$ for all $z_H<z<z_P$. $\mathcal{M}$ is thus an attractor for the (BLS).

A precise study of this sequence of bifurcations is available in \cite{Vidal06, Vidal07}, including the double homoclinic bifurcation that may occur in a certain range of the parameter space. However, here, we focus on the oscillatory behavior of the system in the phase space and control the two types of ejections giving rise to the hysteresis loop: the saddle-node bifurcation and the saddle-node transcritical bifurcation.

\subsection{Global behavior of the system and canard cycles} \label{GlobBehavTT}
 \qquad \newline

The separatrix $\dot{z}=0$ is the plane of equation:
\[
y=\frac{d_2}{a_2-b_2 d_2}
\]
For $0 < d_2 < a_2/b_2$ small enough, this plane separates $\Delta $ near which $\dot{z}<0$ and $\mathcal{M}$ near which $\dot{z}>0$.
\begin{figure}[htb]
\centering
\includegraphics[height=10cm]{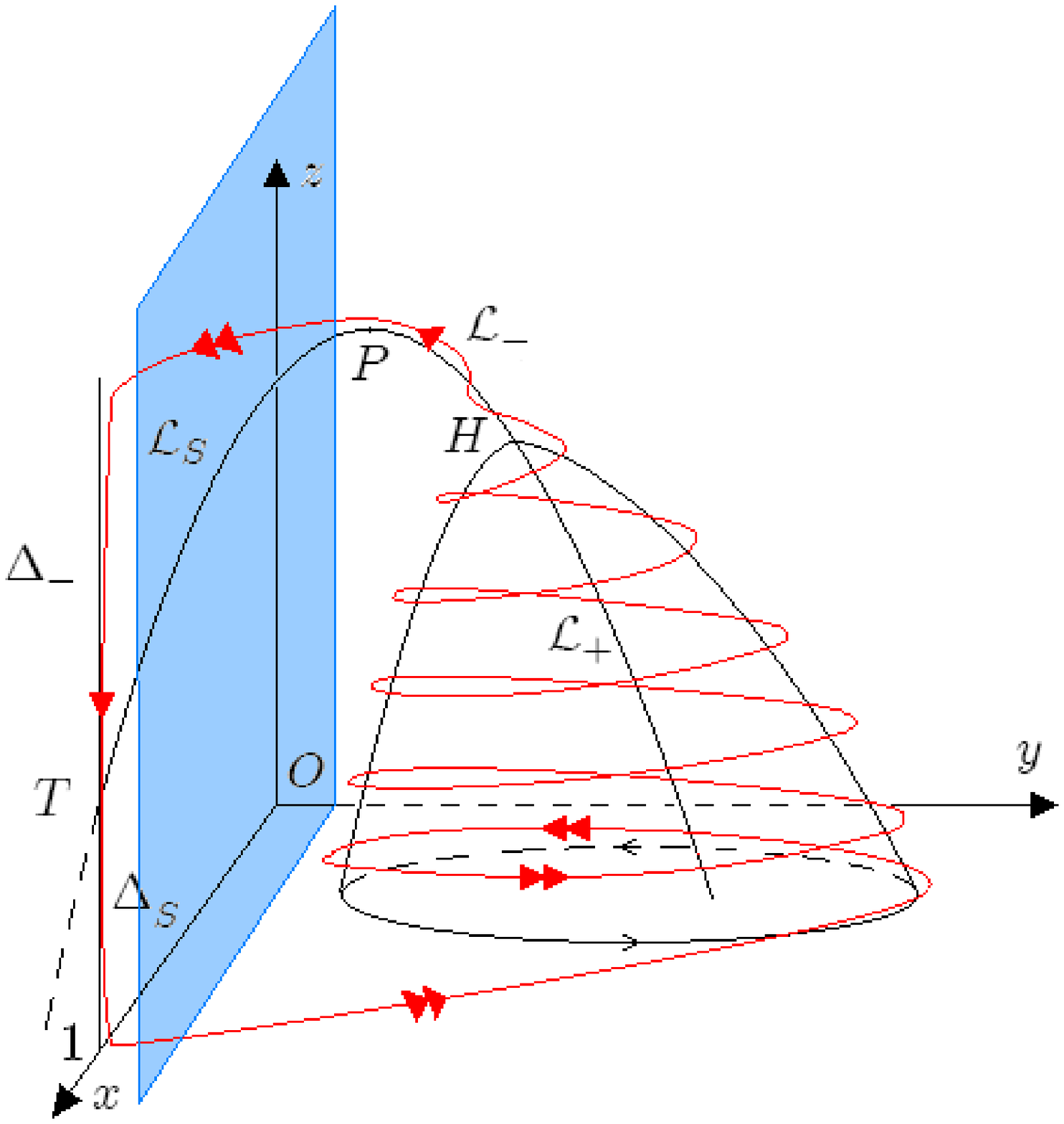}
\caption{Outline of a typical orbit of system \eref{sys2} generating bursting oscillations. Double arrows are added on the fast parts of the orbit and single arrows on the slow parts.}
\label{TTSchemaOrbite}
\end{figure}

Hence, any orbit of system \eref{sys2} starting near $\Delta _-$ goes down along $\Delta _-$. Once $z<z_T$, it undergoes a delay to bifurcation and keeps tracking $\Delta _S$ for a while, although this branch is formed by saddles of the fast dynamics. After this delay, the orbit quickly reaches the vicinity of $\mathcal{M}$. It goes up while spiraling around $\mathcal{M}$, then around $\mathcal{L}_{-}$. As $z$ becomes larger than $z_P$, the orbit reaches the vicinity of $\Delta _-$ and repeats the same sequence of motions. Figure \ref{TTSchemaOrbite} displays the geometric invariants of the fast dynamics and a schematic orbit of \eref{sys2}. In this setting, the delay to the transcritical bifurcation, that the orbit undergoes, qualifies the terminology of canard cycle.

\begin{figure}[htb]
\centering
\includegraphics[height=8.3cm]{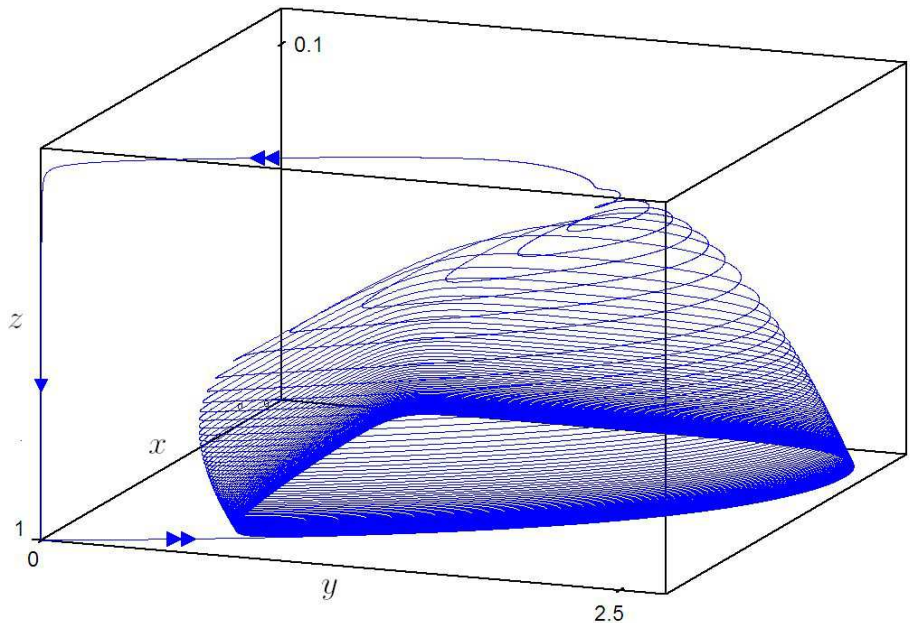}
\caption{Orbit of system \eref{sys2} with the parameter values given in the table of Subsection \ref{GlobBehavTT}. Double arrows are added on the fast parts of the orbit and single arrows on the slow parts.}
\label{TTorbit}
\end{figure}

Figure \ref{TTorbit} displays a numerical simulation of a typical orbit of \eref{sys2} with:
\begin{center}
\begin{tabular}{rclcrclcrclcrcl}
$a_{1}$ & $=$ & $0.8$, & & $b_{1}$ & $=$ & $4$,  & & $d_{1}$ & $=$ & $0.1$,  & & $ \varepsilon $ & $=$ & $0.1$, \\
$a_{2}$ & $=$ & $42$,  && $b_{2}$ & $=$ & $40$, & & $d_{2}$ & $=$ & $1$.    & &                 &     & \\
\end{tabular}
\end{center}

\subsection{Local form near the saddle} \qquad \newline

The most surprising behavior of these orbits is the long time needed to escape a neighbordhood of the point $(1,0,0)$ and the ``exceptionally fast recovery'' of the variable $y$. In this subsection, we perform successively two changes of variables to obtain an appropriate local form near the saddle.

The Jacobian matrix associated with the system \eref{sys1} at the singular point $(x,y,z)=(1,0,0)$ reads:
\begin{equation*}
\begin{pmatrix}
-1 & -\frac{a_{1}}{1+b_{1}}      & 0 \\ 
0  & \frac{a_{1}}{1+b_{1}}-d_{1} & 0 \\ 
0  & 0                           & -\varepsilon d_{2}
\end{pmatrix}
\end{equation*}
It admits $-1$, $\frac{a_{1}}{1+b_{1}}-d_{1} =\frac{G}{1+b_1}$ and $-\varepsilon d_{2}$ as eigenvalues. Under the assumption \eref{StAss}, the singular point $(1,0,0)$ is then a saddle with a two-dimensional stable manifold and a one-dimensional unstable manifold.

Eigenvectors associated with the eigenvalues are:
\begin{center}
\begin{tabular}{|c|c|}
\hline
eigenvalue & eigenvector\\
\hline
$-1$ & $(1,0,0)$ \\
\hline
$\frac{G}{1+b_1}$ & $(-1, \frac{1}{a_1}(G+1+b_1), 0)$ \\
\hline
$-\varepsilon d_{2}$  & $(0,0,1)$ \\
\hline
\end{tabular}
\label{EigVctxyz}
\end{center}
Set:
\[
\alpha=\frac{1}{a_1}(G+1+b_1)=1-\frac{(1+b_1)(d_1-1)}{a_1}
\]
The parameter $\alpha $ is positive under the assumption \eref{StAss} ($G>0$).

Hence, via the change of variables:
\begin{eqnarray}
x &=& 1-X-Y \notag \\
y &=& \alpha Y \notag \\
z &=& Z \notag
\end{eqnarray}
the singular point $(x,y,z)=(1,0,0)$ translates to $(X,Y,Z)=(0,0,0)$. Now, the eigendirections of the saddle coincide with the $(X,Y,Z)$-axis. In this new system of coordinates, \eref{sys2} reads:
\begin{eqnarray}
\dot{X}& =&-(1-X-Y)(X+Y)  \notag \\
& & +Y\left[ -1+\frac{(X-Y)(d_{1}-1)}{1+b_{1}(1+X-Y)}-\frac{a_{2}Z}{1+b_{2}\alpha Y}\right]  \notag \\
\dot{Y}& =&Y\left[ \frac{a_{1}(1+X-Y)}{1+b_{1}(1+X-Y)}-d_{1}-\frac{a_{2}Z}{1+b_{2}\alpha Y}\right]   \label{sys3} \\
\dot{Z}& =&Z\varepsilon \left[ \frac{a_{2}\alpha Y}{1+b_{2}\alpha Y}-d_{2}\right]   \notag 
\end{eqnarray}
This allows to consider $X>0$ in the following.

Actually, an eigenvector associated with each saddle $(0,1)$ of the fast dynamics writes:
\begin{equation}
(\beta (Z), 1)=\left( \frac{a_2 Z (1+b_1)}{1+b_1+G- a_2 Z (1+b_1)} , 1 \right)
\label{VectTT}
\end{equation}
This vector is well-defined in the positive octant as long as:
\begin{equation*}
Z < \frac{1+b_1+G}{a_2 (1+b_1)} = \frac{1}{a_2} \left(1+\frac{G}{1+b_1}\right)
\end{equation*}
Thus, under the assumption \eref{StAss}, it is well-defined at least for $Z\in [0, z_{T}]$ and $\beta(Z)>0$. Note that, for $Z=z_T$, this vector gives the central direction associated with the $0$ eigenvalues of the non hyperbolic point $(0,0)$ of the fast dynamics.

For each $\tilde{Z} \in [0, z_T]$, the vector \eref{VectTT} gives the tangent line to the unstable manifold of the saddle $(0,0)$ of the fast dynamics with $Z=\tilde{Z}$. Thus, we have an approximation for $X$ and $Y$ small of the invariant manifold of \eref{sys3} with $\varepsilon =0$. As this manifold is normally attractive, it persists for $\varepsilon >0$ small enough into a normally attractive manifold of \eref{sys3}. We proceed with another change of variables:
\begin{eqnarray}
u &=& X-\beta (Z) Y \notag \\
y &=& Y \notag \\
z &=& Z \notag
\end{eqnarray}
locally near $X=Y=0, 0 \leq Z \leq z_T$, system \eref{sys3} reads:
\begin{eqnarray}
\dot{u}& = & -u+O(\varepsilon v, u v, u^2)  \notag \\
\dot{v}& = & G' v - a_2 v w + v O(u,v)  \label{sys4} \\
\dot{w}& = & -\varepsilon d_2 w + w O(v^2) \notag 
\end{eqnarray}
where:
\[
G'=\frac{G}{1+b_1}=\frac{a_1}{1+b_1}-d_1
\]
The perturbed attractive manifold in this new set of parameters is approximated by $u=0, 0 \leq Z \leq z_T$.

The flow on the invariant perturbed manifold is conjugated -- locally near $u=v=0$ -- to the system:
\begin{eqnarray}
\dot{v}& =&v \left( \frac{a_1}{1+b_1}-d_1 \right) - a_2 v w \label{sys5} \\
\dot{w}& =&-\varepsilon d_2 w \notag 
\end{eqnarray}
After setting $w=2Y/a_2$ and replacing $\varepsilon d_2$ by $\varepsilon $, we recognize the normal form \eref{DTCBNL} introduced in Subsection \ref{SectDTCBB} with parameter values $G'=(b-1)>0$.

\subsection{Transition function: recovery of the predator} \qquad \newline

Consider the two sections:
\begin{eqnarray}
\Sigma_{in} &=& \left\{(u,v,\delta )|0<v \leq \eta \right\} \\
\Sigma_{out} &=& \left\{(u,\eta,Y)|0<u \leq \xi, 0<w \leq \delta \right\}
\end{eqnarray}
where $\eta, \xi>0$ are small but fixed parameters and $0<\delta < z_T$ (see \eref{zT}) is fixed. 
We note:
\begin{equation*}
U=\left\{(u,v,w )|0<u \leq \xi, 0<v \leq \eta, 0<w \leq \delta \right\} 
\end{equation*}
The preceding analysis shows that, for $\varepsilon$ small enough, an orbit of \eref{sys5} starting from $\Sigma_{in}$ enters $U$ and exits from $U$. Hence, the transition function induced by the flow of \eref{sys5} is well-defined from $\Sigma_{in}$ into $\Sigma_{out}$.

As the typical orbits of system \eref{sys2} enter an exponentially small neighborhood of $\Delta _S$ (see Figures \ref{TTSchemaOrbite}, \ref{TTorbit} and Subsection \ref{GlobBehavTT}), we restrict initial data of \eref{sys5} to be exponentially close to $X=Y=0$ in $\Sigma _{in}$. Thus, consider:
\begin{eqnarray}
u_0=C_1 {\rm e}^{-\frac{k_1}{\varepsilon }} \\
v_0=C_2 {\rm e}^{-\frac{k_1}{\varepsilon }} 
\end{eqnarray}

By direct application of the local analysis made in Subsection \ref{SectDTCBB}, we obtain the leading term of the transition time from $(u_0, v_0, \delta) \in \Sigma_{in}$ to $\Sigma_{out}$ along the flow of \eref{sys5}:
\begin{equation*}
T \underset{\varepsilon \rightarrow 0}{=} \frac{1}{\varepsilon d_2} \left[ \mathcal{W_L} \left( -\frac{a_2 \delta}{G'} {\rm e}^{-\frac{a_2 \delta + k_2 d_2}{G}} \right) + \frac{a_2 \delta + k_2 d_2}{G'} \right] + O(1)
\end{equation*}
The transition function $(u_0,v_0,\delta) \rightarrow (u_{out}(u_0,v_0), \eta, w_{out}(u_0,v_0))$ can be specified as well:
\begin{multline}
u_{out}(C_1 {\rm e}^{-\frac{k_1}{\varepsilon }},C_2 {\rm e}^{-\frac{k_1}{\varepsilon }}) \\
 \underset{\varepsilon \rightarrow 0}{\sim } C_1 {\rm e}^{-\frac{k_1}{\varepsilon}} \left[ -\mathcal{W_L} \left( -\frac{a_2 \delta}{G'} {\rm e}^{-\frac{a_2 \delta + k_2 d_2}{G'}} \right)-\frac{a_2 \delta + k_2 d_2}{d_2 G'} \right]
\end{multline}
\begin{multline}
w_{out}(C_1 {\rm e}^{-\frac{k_1}{\varepsilon }},C_2 {\rm e}^{-\frac{k_1}{\varepsilon }}) \\
 \underset{\varepsilon \rightarrow 0}{= } \delta \exp \left[ -\mathcal{W_L} \left( -\frac{a_2 \delta}{G'} {\rm e}^{-\frac{a_2 \delta + k_2 d_2}{G'}} \right)-\frac{a_2 \delta + k_2 d_2}{d_2 G'} \right]+ O(\varepsilon )
\end{multline}

It is worth noticing that the $u$-component of the transition function tends actually to $0$ as $\varepsilon \rightarrow 0$, but that the $w$-component experiments a basal threshold, even if this threshold is very low.

\subsection{Interpretation in terms of ecological parameters: loss of resilience} \qquad \newline

\vspace{-0.5cm}
In contrast with system \eref{DTCBB}, the two ejection mechanisms in \eref{sys2} are provided by different types of bifurcation of the fast dynamics: a saddle-node transcritical bifurcation at $T$ and a saddle-node bifurcation at $P$. Near $P$, canards could eventually appear for which the orbit tracks the manifold $\mathcal{L}_S$ formed by saddle of the fast dynamics. However, this type of canards can occur only if $P$ is near the plane $\dot{z}=0$ (within a $O(\sqrt{\varepsilon})$ distance). But the assumption that this plane separates properly $\Delta $ and the manifold $\mathcal{M}$ ($d_2$ small enough) generally forbids such behavior.

Hence, any orbit in the strictly positive octant displays, sooner or later, the motion described in Subsection \ref{GlobBehavTT} as there exists a globally attractive periodic orbit of this type (see \cite{Vidal06, Vidal07}). Thus, the entry of this periodic orbit -- as the entry of any orbit after spiraling around $\mathcal{M}$ -- in a given neighborhood of $\Delta $ occurs for $z=z_P+O(\varepsilon ^{2/3})$. The delay to bifurcation after the passage in the vicinity of $T$ is then mostly characterized by $z_P - z_T$. In fact, as calculated for \eref{DTCBB} in \eref{paramk}, the contraction exponent between $\Delta $ and the orbit is chiefly determined by $z_P - z_T$.

The ``exceptionally fast recovery'' that we have pointed out therein may have some connection with the notion of ``loss of resilience'' which has been discussed in several dynamical models related to Ecology \cite{H, LWH, Martin}. The long return times associated with a loss of resilience are caused by slow dynamics near the unstable equilibrium. We expect that new developments of this subject could possibly benefit of the notion of canards and canard cycles.

\vspace{0.5cm}

\textbf{Acknowledgements:} The numerical simulations were performed using XPP-AUT:
\href{http://www.math.pitt.edu/\~{}bard/xpp/xpp.html}{\tt http://www.math.pitt.edu/\~{}bard/xpp/xpp.html}

\end{document}